# Confluent hypergeometric expansions of the solutions of the double-confluent Heun equation

T.A. Ishkhanyan[1,2], V.A. Manukyan[2], A.H. Harutyunyan[2], A.M. Ishkhanyan[2,3]

[1]*Laboratoire Interdisciplinaire Carnot de Bourgogne, Université de Bourgogne Franche-Comté, 9 Av. Alain Savary, 21078 Dijon, France*
[2]*Institute for Physical Research, NAS of Armenia, Ashtarak-2, 0203 Ashtarak, Armenia*
[3]*Russian-Armenian University, H. Emin 123, 0051 Yerevan, Armenia*
E-mail: tishkhanyan@gmail.com

**Abstract.** Several expansions of the solutions of the double-confluent Heun equation in terms of the Kummer confluent hypergeometric functions are presented. Three different sets of these functions are examined. Discussing the expansions without a pre-factor, it is shown that two of these functions provide expansions the coefficients of which obey three-term recurrence relations, while for the third confluent hypergeometric function the corresponding recurrence relation generally involves five-terms. The latter relation is reduced to a three-term one only in the case when the double-confluent Heun equation degenerates to the confluent hypergeometric equation. The conditions for obtaining finite sum solutions via termination of the expansions are discussed. The possibility of constructing expansions of different structure using certain equations related to the double-confluent Heun equation is discussed. An example of such expansion derived using the equation obeyed by a function involving the derivative of a solution of the double-confluent Heun equation is presented. In this way, expansions governed by three- or more term recurrence relations for expansion coefficients can be constructed. An expansion with coefficients obeying a seven-term recurrence relation is presented. This relation is reduced to a five-term one if the additional singularity of the equation obeyed by the considered auxiliary function coincides with a singularity of the double-confluent Heun equation.

## 1. Introduction

The double-confluent Heun equation is one of the four confluent reductions of the general Heun equation [1]. This is a linear second order ordinary differential equation having two irregular singularities of rank 1, each formed by coalescence of a pair of regular singularities of the general Heun equation [2,3]. If the singularities are placed at $z = 0$ and $\infty$, the double-confluent Heun equation is written as

$$\frac{d^2u}{dz^2} + \left(\frac{\delta}{z^2} + \frac{\gamma}{z} + \varepsilon\right)\frac{du}{dz} + \frac{\alpha z - q}{z^2}u = 0. \qquad (1)$$

Note that this form of the equation slightly differs from that adopted in [4], namely, here we have introduced the parameter $\varepsilon$ instead of the unity used there. This notation is useful in order to distinguish the Whittaker-Ince limit [5] of the double-confluent Heun equation as a particular case achieved within this form by the simple choice $\varepsilon = 0$. We will see that the notation is also helpful to explicitly see that the presented expansions do not apply to this limiting case. Several expansions in terms of regular and irregular confluent hypergeometric functions and Bessel functions which apply to this particular limit are discussed in [6,7].

Expansions of the solutions of the general Heun equation in terms of special functions initiated by Svartholm and Erdélyi [8-9] have been extended to more general equations [10] as well as to the confluent forms of the Heun equation by many authors (see, e.g., [2,3] and references therein). Different special functions have been used as expansion functions. In the context of the double-confluent Heun equation, the regular and irregular confluent hypergeometric functions [2,3,6,7,11], the Gauss hypergeometric functions [2,3], Bessel functions [6,7], Coulomb wave functions [12,13], incomplete beta and gamma functions [14,15], and other special functions have been applied. Using the properties of the derivatives of the solutions of the Heun equations [16-20], expansions in terms of higher transcendental functions such as the Appell and Goursat generalized hypergeometric functions are also possible [17,19-20]. The expansions in terms of different special functions may have different regions of convergence, hence, may be applicable to diverse physical situations [2,3].

In the present paper, we discuss the expansions of the solutions of the double-confluent Heun equation in terms of the Kummer confluent hypergeometric functions (confluent hypergeometric functions of the first kind). The presented expansions are of general interest both for mathematics and physics, and may have useful applications in many branches of physics varying from quantum mechanics, atomic, molecular and optical physics to particle physics or astrophysics and cosmology [21].

Starting from the recurrence relations that these functions obey for chosen forms of the dependence on the summation variable, we first present the direct expansions derived by means of substitution of the expansion ansatz without a pre-factor into the double-confluent Heun equation. The coefficients of such expansions for two of applied hypergeometric functions generally obey three-term recurrence relations. For one of the discussed forms of the involved confluent hypergeometric functions a two-term recurrence relation is also possible for a particular choice of the involved parameters. In this case the coefficients of the expansion are explicitly calculated. Another expansion governed by a three-term recurrence relation for expansion coefficients, involving a third type of hypergeometric functions, is applicable only in the degenerate case when the double-confluent Heun equation is reduced to the confluent hypergeometric equation by a simple transformation of the dependent variable. However, we show that it is possible to construct, with the same expansion function, a different expansion applicable in the general non-degenerate case. The coefficients of the latter expansion, however, obey a five-term recurrence relation. The conditions for termination of the presented expansions are discussed.

In the last section we note that different expansions can be constructed by employment of other equations related to the double-confluent Heun equation. As an example, we consider a second-order differential equation obeyed by a function involving the derivative of a solution of the double-confluent Heun equation. Using the expansions of the solutions of this equation in terms of the confluent hypergeometric functions, one can construct confluent hypergeometric expansions of the solutions of the double-confluent Heun equation, which involve pre-factors consisting of elementary functions. We present an example of such expansion which generally is governed by a seven-term recurrence relation for the expansion coefficients. Two particular cases when the relation involves five terms are mentioned. This happens when the additional singularity of the equation obeyed by the considered function involving the derivative of a solution of the double-confluent Heun equation coincides with already existing singularity of the double-confluent Heun equation.

## 2. Direct expansions through the double-confluent Heun equation

Consider the direct expansions of the solutions of the double-confluent Heun equation

$$u = \sum_n a_n u_n, \quad u_n = {}_1F_1(\alpha_n; \gamma_n; s_0 z), \qquad (2)$$

where the Kummer confluent hypergeometric functions $u_n = {}_1F_1(\alpha_n; \gamma_n; s_0 z)$ have one of the following forms: ${}_1F_1(\alpha_0 + n; \gamma_0 + n; s_0 z)$, ${}_1F_1(\alpha_0 + n; \gamma_0; s_0 z)$ and ${}_1F_1(\alpha_0; \gamma_0 + n; s_0 z)$. We will see that the functions of the second type provide an expansion governed by a three-term recurrence relation for expansion coefficients only in the degenerate case $\delta = 0$ when the double-confluent Heun equation is reduced to the confluent hypergeometric equation by a simple transformation of the dependent variable.

The chosen expansion functions obey the confluent hypergeometric equation

$$u_n'' + \left(\frac{\gamma_n}{z} - s_0\right)u_n' - \frac{\alpha_n s_0}{z}u_n = 0, \tag{3}$$

so that the substitution of Eqs. (2) and (3) into Eq. (1) yields

$$\sum_n a_n \left[\left(\frac{\delta}{z^2} + \frac{\gamma - \gamma_n}{z} + \varepsilon + s_0\right)u_n' + \frac{(\alpha + \alpha_n s_0)z - q}{z^2}u_n\right] = 0 \tag{4}$$

or $\quad \sum_n a_n \left[\left((\varepsilon + s_0)z^2 + (\gamma - \gamma_n)z + \delta\right)u_n' + \left((\alpha + \alpha_n s_0)z - q\right)u_n\right] = 0. \tag{5}$

The discussed confluent hypergeometric functions obey different recurrence relations depending on the particular form of the dependence on the summation variable $n$. Consider the three possible cases separately.

**2.1.** Let $\alpha_n = \alpha_0 + n$, $\gamma_n = \gamma_0 + n \Rightarrow u_n = {}_1F_1(\alpha_0 + n; \gamma_0 + n; s_0 z)$.

The differentiation rule for the Kummer confluent hypergeometric functions suggests a starting recurrence relation: $u_n' = s_0(\alpha_n / \gamma_n)u_{n+1}$. Since the terms $z^2 u_n'$, $z u_n'$ and $z u_n$ are not expressed as linear combinations of functions $u_n$ one may first try to cancel these terms. However, we then get that should be $s_0 = -\varepsilon$ and $\alpha_n = \alpha/\varepsilon$, $\gamma_n = \gamma$ for any $n$, i.e. $\alpha_n$, $\gamma_n$ should not depend on $n$. Hence, we conclude that in this way no expansion is constructed.

Another possibility opens if we put $s_0 = -\varepsilon$ (thereby canceling the term proportional to $z^2$) and demand

$$(\gamma - \gamma_n)z + \delta = A f(n+1, z), \quad (\alpha + \alpha_n s_0)z - q = f(n, z), \tag{6}$$

where $A$ is a constant and $f(n, z)$ is a liner function of $n$. This is possible only if $A = 1/\varepsilon$, $\gamma_0 = 1 + \alpha_0 + \gamma - \alpha/\varepsilon$ and $q = -\delta \varepsilon$. Eq. (5) then becomes

$$\sum_n a_n \left(-f(n+1, z)\frac{\alpha_n}{\gamma_n}u_{n+1} + f(n, z)u_n\right) = 0, \tag{7}$$

from where we get the following simple two-term recurrence relation:

$$a_n - \frac{\alpha_{n-1}}{\gamma_{n-1}} a_{n-1} = 0, \tag{8}$$

Thus, for $q = -\delta\varepsilon$ applies the expansion

$$u = \sum_n a_n \cdot {}_1F_1(\alpha_0 + n; \gamma_0 + n; -\varepsilon z), \tag{9}$$

with $\quad a_n = \frac{(\alpha_0)_n}{(\gamma_0)_n}, \quad \gamma_0 = 1 + \alpha_0 + \gamma - \alpha/\varepsilon, \tag{10}$

where $(...)_n$ is the Pochhammer symbol. In general, this is a double-sided infinite series which is valid if $\varepsilon \neq 0$ and $\gamma_0$ is not zero or a negative integer. The series is right-hand side terminated if $\alpha_0 = 0$ or $\alpha_0 = -N$ for some positive integer $N$. Choosing $\alpha_0 = 0$ and changing $n \to -n$ we arrive at the following polynomial expansion for $q = -\delta\varepsilon$:

$$u = \sum_{n=0}^{\infty} \frac{(1-\gamma_0)_n}{n!} \cdot {}_1F_1(-n; \gamma_0 - n; -\varepsilon z), \quad \gamma_0 = 1 + \gamma - \alpha/\varepsilon, \tag{11}$$

which is a particular case of the known hypergeometric expansions [2-3].

The additional restriction $q = -\delta\varepsilon$ imposed on the parameters of the double-confluent Heun equation (1) is avoided if we use the recurrence relation

$$z(u_n' - s_0 u_n) = (\gamma_n - 1)(u_{n-1} - u_n). \tag{12}$$

Indeed, we again put $s_0 = -\varepsilon$ and demand [compare with Eq. (6)]

$$\alpha + \alpha_n s_0 = -s_0 [(\gamma - \gamma_n)]. \tag{13}$$

This equation is fulfilled if

$$\alpha_0 = \gamma_0 - \gamma + \alpha/\varepsilon. \tag{14}$$

With this, Eq. (5) is rewritten as

$$\sum_n a_n \left[ (\gamma - \gamma_n) z(u'_n - s_0 u_n) + \delta u'_n - q u_n \right] = 0. \tag{15}$$

Substituting Eq. (12), we get

$$\sum_n a_n \left[ (\gamma - \gamma_n)(\gamma_n - 1)(u_{n-1} - u_n) + \delta s_0 \frac{\alpha_n}{\gamma_n} u_{n+1} - q u_n \right] = 0, \tag{16}$$

from where we derive the following three-term recurrence relation for the coefficients of the expansion:

$$R_n a_n + Q_{n-1} a_{n-1} + P_{n-2} a_{n-2} = 0, \tag{17}$$

where

$$R_n = (\gamma - \gamma_n)(\gamma_n - 1), \quad Q_n = -R_n - q, \quad P_n = -\varepsilon \delta \alpha_n / \gamma_n. \tag{18}$$

For left-hand side termination of this series at $n = 0$ should be $R_0 = 0$. This is the case if

$$\gamma_0 = \gamma \quad (\Rightarrow \alpha_0 = \alpha/\varepsilon) \tag{19}$$

($\gamma_0 = 1$ is forbidden because of division by zero: $P_{-1} \sim 1/\gamma_{-1}$, $\gamma_{-1} = 0$). Thus, finally, the expansion is explicitly written as

$$u = \sum_{n=0}^{\infty} a_n \cdot {}_1F_1((\alpha/\varepsilon) + n; \gamma + n; -\varepsilon z), \tag{20}$$

and the coefficients of the recurrence relation (17) are explicitly given as

$$R_n = -n(\gamma + n - 1), \quad Q_n = n(\gamma + n - 1) - q, \quad P_n = -\delta \frac{\alpha + \varepsilon n}{\gamma + n}. \tag{21}$$

This expansion has been derived in [6] (Eqs. (57a) and (57b)). Like the expansion (9)-(10), it is applicable if $\varepsilon \neq 0$ and $\gamma$ is not zero or a negative integer. Note that $P_n$ is identically zero if $\delta = 0$, so that in this case the recurrence relation becomes two-term, and the coefficients $a_n$ are explicitly calculated in terms of gamma-functions. However, the case $\delta = 0$ is degenerate since in this case the singularity at $z = 0$ becomes regular and the double-confluent Heun equation (1) is readily reduced to the confluent hypergeometric equation.

The series is terminated from the right-hand side if $a_{N+1} = 0$ and $a_{N+2} = 0$ for some non-negative integer $N$. The first condition is fulfilled when $Q_N a_N + P_{N-1} a_{N-1} = 0$, while for a non-zero $a_N$ the second condition results in the equation $P_N = 0$, which is satisfied, apart from the degenerate case $\delta = 0$, if

$$\alpha/\varepsilon = -N. \tag{22}$$

The equation $Q_N a_N + P_{N-1} a_{N-1} = 0$ (or, equivalently, $a_{N+1} = 0$) presents a polynomial equation of the degree $N + 1$ for the accessory parameter $q$, defining, in general, $N + 1$ values of $q$ for which the termination of the series occurs. Note, finally, that the resulting finite-sum solution is a polynomial in $z$.

**2.2.** Let $\alpha_n = \alpha_0 + n$, $\gamma_n = \gamma_0 = \text{const} \Rightarrow u_n = {}_1F_1(\alpha_0 + n; \gamma_0; s_0 z)$.

In this case we have the following recurrence relations

$$z u'_n = \alpha_n (u_{n+1} - u_n), \quad s_0 z u_n = (\alpha_n - \gamma_0) u_{n-1} + (\gamma_0 - 2\alpha_n) u_n + \alpha_n u_{n+1}. \tag{23}$$

Combining these equations we also have

$$s_0 z^2 u'_n = \alpha_n (s_0 z u_{n+1} - s_0 z u_n) = \\ \alpha_n \left( (\alpha_n + 1) u_{n+2} + (\gamma_0 - 3\alpha_n - 2) u_{n+1} + (3\alpha_n - 2\gamma_0 + 1) u_n + (\gamma_0 - \alpha_n) u_{n-1} \right). \tag{24}$$

However, the derivative $u'_n$ alone is not expressed as a linear combination of functions $u_n$. Hence, with this form of expansion functions we are able to treat only the case $\delta = 0$ (see Eq. (5)). This is, however, a degenerate case since in this case the double-confluent Heun equation (1) is known to reduce to the confluent hypergeometric equation by the transformation of the dependent variable $u = z^s v(z)$. Accordingly, the general solution of the double-confluent Heun equation in this case is written as

$$u = z^s \left( C_1 \cdot {}_1F_1(s + \alpha/\varepsilon; \gamma_0; -\varepsilon z) + C_2 U(s + \alpha/\varepsilon; \gamma_0; -\varepsilon z) \right), \tag{25}$$

where $U$ is the Tricomi confluent hypergeometric function, $C_{1,2} = \text{const}$, and

$$s = (\gamma_0 - \gamma)/2, \quad \gamma_0 = 1 + \sqrt{(1-\gamma)^2 + 4q} \tag{26}$$

Though the case $\delta = 0$ is degenerate, for completeness of the treatment and for reference purposes let us present the corresponding expansion. Substituting Eqs. (23) and (24) into Eq. (5) we obtain a four-term recurrence relation for coefficients $a_n$. However, if we put $s_0 = -\varepsilon$, thereby removing the $z^2$-term in the coefficient of $u'_n$ in Eq. (5) the four-term recurrence relation becomes a three-term one:

$$R_n a_n + Q_{n-1} a_{n-1} + P_{n-2} a_{n-2} = 0, \tag{27}$$

$$R_n = (\alpha_n - \gamma_0)(\alpha_n - \alpha/\varepsilon), \tag{28}$$

$$Q_n = (\alpha_n - \alpha/\varepsilon)(\gamma_0 - 2\alpha_n) - \alpha_n(\gamma - \gamma_0) - q, \tag{29}$$

$$P_n = \alpha_n (\alpha_n + \gamma - \gamma_0 - \alpha/\varepsilon). \tag{30}$$

The initial conditions for left-hand side termination of the derived series at $n = 0$ are $a_{-2} = a_{-1} = 0$. Hence, for that should be $R_0 = 0$. This is the case if $\alpha_0 = \alpha/\varepsilon$ or $\alpha_0 = \gamma_0$.

Thus, the expansion

$$u = \sum_{n=0}^{\infty} a_n \cdot {}_1F_1(\alpha_0 + n; \gamma_0; -\varepsilon z) \tag{31}$$

applies only if $\delta = 0$. The coefficients of the expansion are defined by relations (27)-(30), where $\gamma_0$ is an arbitrary constant ($\gamma_0 \neq -n$) and $\alpha_0 = \alpha/\varepsilon$ or $\alpha_0 = \gamma_0$. The series is right-side terminated at some $n = N$ if $a_N \neq 0$ and $a_{N+1} = a_{N+2} = 0$. Then, should be $P_N = 0$. If $\alpha_0 = \alpha/\varepsilon$, this condition is satisfied when

$$\alpha/\varepsilon = -N \quad \text{or} \quad \gamma - \gamma_0 = -N. \tag{32}$$

If $\alpha_0 = \gamma_0$, the only possibility, since $\gamma_0$ should not be a negative integer number, is

$$\gamma - \alpha/\varepsilon = -N. \tag{33}$$

Again, for each of these cases there exist $N+1$ values of $q$ for which the termination occurs. These values are determined from the equation $a_{N+1} = 0$ (or $Q_N a_N + P_{N-1} a_{N-1} = 0$).

However, the above development can be extended to include the case $\delta \neq 0$. This is achieved by employing the first recurrence relation of Eqs. (23) after Eq. (5) is multiplied by $z$. In this way we get a five-term recurrence relation:

$$T_n a_n + S_{n-1} a_{n-1} + R_{n-2} a_{n-2} + Q_{n-3} a_{n-3} + P_{n-4} a_{n-4} = 0, \tag{34}$$

with rather cumbersome coefficients. The coefficient at $a_n$ is given as

$$T_n = (\alpha_n - \gamma_0)(\alpha_n - 1 - \gamma_0)(\alpha - \varepsilon\alpha_n), \tag{35}$$

so that the series (2) may left-hand side terminate at $n=0$ ($T_0 = 0$) if $\alpha_0 = \gamma_0, 1+\gamma_0, \alpha/\varepsilon$. However, only the last two cases are permissible because the first choice, $\alpha_0 = \gamma_0$, produces division by zero when calculating $a_1$. At $\alpha_0 = 1+\gamma_0$ the expansion functions are reducible confluent hypergeometric functions, hence, the most interesting is the choice

$$\alpha_0 = \alpha/\varepsilon. \tag{36}$$

With this, the expansion reads

$$u = \sum_{n=0}^{\infty} a_n \cdot {}_1F_1(\alpha/\varepsilon + n; \gamma_0; -\varepsilon z) \tag{37}$$

and the coefficients of the five-term recurrence relation (34) are explicitly written as

$$T_n = -n(\alpha + (n-1-\gamma_0)\varepsilon)(\alpha + (n-\gamma_0)\varepsilon), \tag{38}$$

$$S_n = \left(\alpha(\gamma - \gamma_0 + 4n) + \varepsilon(q + n(\gamma - 3\gamma_0 + 4n - 2))\right)(\alpha + (n-\gamma_0)\varepsilon), \tag{39}$$

$$R_n = -(T_n + S_n + Q_n + P_n), \tag{40}$$

$$Q_n = (\alpha + n\varepsilon)\left(\alpha(3(\gamma - \gamma_0) + 4n) + \varepsilon\left((3\gamma + 4n + 2)n - (\gamma + 5n + 2)\gamma_0 + \gamma_0^2 + q + 2\gamma + \delta\varepsilon\right)\right), \tag{41}$$

$$P_n = -(n + \gamma - \gamma_0)(\alpha + \varepsilon n)(\alpha + \varepsilon(1+n)). \tag{42}$$

Note that here $\gamma_0$ is a free parameter which can be employed to simplify the coefficients. Right-hand side termination is possible for some $N = 0, 1, 2, ...$ if $a_N \neq 0$, $P_N = 0$, that is

$$\gamma_0 = \gamma + N \quad \text{or} \quad \alpha = -\varepsilon N \quad \text{or} \quad \alpha = -\varepsilon(N+1), \tag{43}$$

and $a_{N+1} = a_{N+2} = a_{N+3} = 0$. The last equations impose additional restrictions on the parameters of the double-confluent Heun equation. We note that in all cases the termination of the series results in polynomial solutions.

**2.3.** Let $\alpha_n = \alpha_0 = \text{const} \neq 0$, $\gamma_n = \gamma_0 + n \Rightarrow u_n = {}_1F_1(\alpha_0; \gamma_0 + n; s_0 z)$.

For these functions the following recurrence relations are known:

$$zu'_n = (\gamma_n - 1)(u_{n-1} - u_n), \quad u'_n = s_0\left(u_n - \left(1 - \frac{\alpha_0}{\gamma_n}\right)u_{n+1}\right). \tag{44}$$

Since in this case $z^2 u'_n$, and $zu_n$ are not expressed as linear combinations of functions $u_n$ we demand the coefficients of these terms in Eq. (5) to be equal to zero. It is seen that should be $s_0 = -\varepsilon$ and $\alpha_0 = \alpha/\varepsilon$, so that Eq. (5) is rewritten as

$$\sum_n a_n\left[((\gamma - \gamma_n)z + \delta)u'_n - qu_n\right] = 0. \tag{45}$$

Substituting the relations (34) into this equation we get the following recurrence relation:

$$R_n a_n + Q_{n-1} a_{n-1} + P_{n-2} a_{n-2} = 0, \tag{46}$$

where $\quad R_n = (\gamma - \gamma_n)(\gamma_n - 1), \quad Q_n = -R_n - \delta\varepsilon - q, \quad P_n = \delta(\varepsilon - \alpha/\gamma_n). \tag{47}$

For left-side termination of the derived series at $n=0$ should be $R_0 = 0$. This is the case only if $\gamma_0 = \gamma$ ($\gamma_0 = 1$ is forbidden because of division by zero: $P_1 \sim 1/\gamma_{-1}$, $\gamma_{-1} = 0$). Thus, the expansion is finally written as

$$u = \sum_{n=0}^{\infty} a_n \cdot {}_1F_1(\alpha/\varepsilon; \gamma + n; -\varepsilon z) \tag{48}$$

and the coefficients of the recurrence relation (46) are explicitly given as (compare with (21))

$$R_n = -n(\gamma + n - 1), \quad Q_n = n(\gamma + n - 1) - \varepsilon\delta - q, \quad P_n = \delta\left(\varepsilon - \frac{\alpha}{\gamma + n}\right). \tag{49}$$

This expansion is equivalent to those derived earlier in [6,7]. It is applicable if $\alpha \neq 0$, $\varepsilon \neq 0$ and $\gamma$ is not zero or a negative integer. The series is right-hand side terminated for some non-negative integer $N$ if $P_N = 0$ or, apart from the degenerate case $\delta = 0$, if

$$\gamma - \alpha/\varepsilon = -N. \tag{50}$$

The termination occurs for $N+1$ values of the accessory parameter $q$ defined from the equation $a_{N+1} = 0$. Using Kummer's identity $_1F_1(a;b;-z) = e^{-z}\,_1F_1(b-a;b;z)$, we note that the resultant finite-sum solutions are quasi-polynomials of exponent $-\varepsilon$ and degree $N$, namely, a product of $e^{-\varepsilon z}$ and a polynomial in $z$ of the degree $N$.

### 3. Other expansions

Several other expansions can be constructed. For instance, the way we have derived the expansion (37)-(42), ruled by a five-term recurrence relation for expansion coefficients, suggests that one can construct expansions in terms of the functions $u_n = \,_1F_1(\alpha_0 + n; \gamma_0; s_0 z)$ if the starting double-confluent Heun equation is preliminary changed by a transformation of the form $u = \phi(x)w(z(x))$ with rational $\phi(x)$ and $z(x)$. Other expansions of different structure can be derived starting from other equations related to the double-confluent Heun equation, e.g., the third order equation obeyed by the square of a solution of the double-confluent Heun equation or the equation obeyed by a function involving the derivative of a solution of the double-confluent Heun equation.

To demonstrate the approach, consider the equation obeyed by the function $v = z^\gamma e^{\varepsilon z - \delta/z} u'$:

$$v_{zz} - \left(\frac{\delta}{z^2} + \frac{\gamma - 2}{z} + \varepsilon + \frac{1}{z - q/\alpha}\right)v_z + \frac{\alpha z - q}{z^2}v = 0. \tag{51}$$

It is seen that this equation has an additional regular singularity located at the point $z_0 = q/\alpha$. This is an apparent singularity the emergence of which when considering the equations obeyed by functions involving the derivative of a solution of the five Heun equations has been noticed in several occasions (see, e.g., [22-26]).

It is readily understood that the expansion of the solution of Eq. (51) in terms of the confluent hypergeometric functions:

$$v = \sum_n a_n v_n, \quad v_n = \,_1F_1(\alpha_n; \gamma_n; s_0 z), \tag{52}$$

results in a corresponding expansion of the solution of the double-confluent Heun equation. Indeed, from Eq. (1) we have

$$u = -\frac{z^2}{\alpha z - q}\left[\frac{u''}{u'} + \left(\frac{\delta}{z^2} + \frac{\gamma}{z} + \varepsilon\right)\right]u'. \tag{53}$$

Hence, since $u' = v \cdot z^{-\gamma} e^{-\varepsilon z + \delta/z}$,

$$u = -\frac{z^2}{\alpha z - q}\left[\frac{u''}{u'} + \left(\frac{\delta}{z^2} + \frac{\gamma}{z} + \varepsilon\right)\right]z^{-\gamma} e^{-\varepsilon z + \delta/z} v(z) \tag{54}$$

or, after simplification,

$$u = -\frac{z^{2-\gamma} e^{-\varepsilon z + \delta/z}}{\alpha z - q} v'(z). \tag{55}$$

Thus, the series (52) results in the expansion

$$u = -\frac{z^{2-\gamma}e^{-\varepsilon z+\delta/z}}{\alpha z - q}\sum_n a_n \frac{s_0 \alpha_n}{\gamma_n}\,{}_1F_1(\alpha_n+1;\gamma_n+1;s_0 z). \tag{56}$$

An immediate observation now is that in two particular cases, namely, when $q=0$ or $\alpha=0$, the additional singularity of Eq. (51) coincides with the singularities of Eq. (1) (with $z=0$ or $z=\infty$, respectively), so that in these two cases Eq. (51) is a double-confluent Heun equation with altered parameters. Hence, one can directly use the above expansions **2.1**, **2.3** to construct series governed by three-term recurrence relations for expansion coefficients or to employ the expansion of section **2.2** to derive a series with a five-term recurrence relation for the coefficients. As regards the general case of non-zero $\alpha, q$, the series (52) is constructed essentially in the same manner as the above expansions **2.1**, **2.2**, **2.3**. However, this time the recurrence relation for the successive coefficients of the expansion involves more terms. For example, for the expansion functions $v_n = {}_1F_1(\alpha_0+n;\gamma_0;s_0 z)$ one multiplies Eq. (51) by $z$ and employs the recurrence relations (23) to get an expansion with a seven-term recurrence relation for the coefficients ($s_0 = \varepsilon$):

$$A_n a_n + B_{n-1} a_{n-1} + T_{n-2} a_{n-2} + S_{n-3} a_{n-3} + R_{n-4} a_{n-4} + Q_{n-5} a_{n-5} + P_{n-6} a_{n-6} = 0. \tag{57}$$

The coefficients $A_n, B_n, T_n, S_n, R_n, Q_n, P_n$ of this relation are rather cumbersome. Relatively simple are the first and the last ones, which read

$$A_n = (\alpha_n - \gamma_0)(\alpha_n - \gamma_0 - 1)(\alpha_n - \gamma_0 - 2)(\alpha_n + \alpha/\varepsilon), \tag{58}$$

$$P_n = \alpha_n(\alpha_n+1)(\alpha_n+2)(\alpha_n+1-\gamma_0-\delta+\alpha/\varepsilon). \tag{59}$$

The series is left-hand side terminated at $n=0$ if $A_0 = 0$, which gives only 2 possible values for $\alpha_0$: $\alpha_0 = \gamma_0 + 2$ and $\alpha_0 = -\alpha/\varepsilon$ (the cases $\alpha_0 = \gamma_0$ and $\alpha_0 = \gamma_0 + 1$ are excluded because they produce divisions by zero in calculating $a_1, a_2$). The conditions for right-hand side termination define an over-determined system of equations that can be satisfied only occasionally. Indeed, the right-hand side termination is possible for some $N = 0,1,2...$ if $P_N = 0$ (that is $\alpha_0 + N = 0, -1, -2, \gamma_0 + \delta - \alpha/\varepsilon - 1$) and, in addition, if six more restrictions hold for the five remaining free parameters: $a_{N+1} = a_{N+2} = a_{N+3} = a_{N+4} = a_{N+5} = a_{N+6} = 0$.

## 4. Summary

Thus, we have presented several expansions of the solutions of the double-confluent Heun equation in terms of the Kummer confluent hypergeometric functions. The expansions apply if $\varepsilon \neq 0$ and $\gamma$ is not zero or a negative integer (some expansions are applicable if additionally $\alpha \neq 0$). We have examined three different sets of the Kummer functions.

In section **2** we have presented the expansions for which the coefficients generally obey three-term recurrence relations. We have indicated a particular case for which a two-term relation is derived. Furthermore, we have seen that for one of the discussed confluent hypergeometric functions ($u_n = {}_1F_1(\alpha_0+n;\gamma_0;s_0z)$) the recurrence relation generally involves five-terms. This relation is reduced to a three-term one only in the case when the double-confluent Heun equation degenerates to the confluent hypergeometric equation. We have discussed the conditions for obtaining finite sum solutions via termination of the series.

Comparing the expansions presented in section **2** with similar expansions (in terms of the Kummer functions of the same form as the ones used here) of the solutions of the single-confluent Heun equation [27], we see that the main difference consists in that in the case of the double-confluent Heun equation the expansion in terms of functions ${}_1F_1(\alpha_0+n;\gamma_0;s_0z)$ with a three-term recurrence relation for the coefficients applies only in the degenerate case $\delta = 0$. Another difference for the series ruled by three-term recurrence relations is noted in the conditions for the right-hand side termination. In the case of the double-confluent Heun equation these conditions are less flexible than in the case of the confluent Heun equation.

Finally, we have noted that other expansions having different structure can be constructed using certain equations related to the double-confluent Heun equation. We have presented an example of such an expansion applying the equation obeyed by a function involving the derivative of a solution of the double-confluent Heun equation. The coefficients of this expansion in general obey a seven-term recurrence relation. This relation is reduced to a five-term one if the additional singularity of the equation obeyed by the considered function involving the derivative of a solution of the double-confluent Heun equation coincides with a singularity of the double-confluent Heun equation.

**Acknowledgments**

This research has been conducted within the scope of the International Associated Laboratory IRMAS (CNRS-France & SCS-Armenia). The work has been supported by the Armenian State Committee of Science (SCS Grant No. 15T-1C323) and by the Armenian National Science and Education Fund (ANSEF Grant No. PS-4558). T.A. Ishkhanyan acknowledges the support from SPIE through a 2017 Optics and Photonics Education Scholarship and thanks the French Embassy in Armenia for a doctoral grant.

**References**
[1] K. Heun, Math. Ann. **33** (1889) 161.
[2] A. Ronveaux (ed.), Heun's Differential Equations (Oxford Univ. Press, London, 1995).
[3] S.Yu. Slavyanov and W. Lay, Special functions (Oxford Univ. Press, Oxford, 2000).
[4] F.W.J. Olver, D.W. Lozier, R.F. Boisvert, and C.W. Clark (eds.), *NIST Handbook of Mathematical Functions* (Cambridge Univ. Press, New York, 2010).
[5] E.L. Ince, Proc. London Math. Soc. **23** (1925) 56.
[6] B.D.B. Figueiredo, J. Math. Phys. **48** (2007) 013503.
[7] L.J. El-Jaick and B.D.B. Figueiredo, J. Math. Phys. **49** (2008) 083508.
[8] N. Svartholm, Math. Ann. **116** (1939) 413.
[9] A. Erdélyi, Q. J. Math. (Oxford), **15** (1944) 62.
[10] D. Schmidt, J. Reine Angew. Math. **309** (1979) 127.
[11] Th. Kurth and D. Schmidt, SIAM J. Math. Anal. **17** (1986) 1086.
[12] L.J. El-Jaick and B.D.B. Figueiredo, J. Phys. A **46** (2013) 085203.
[13] E.W. Leaver, J. Math. Phys. **27** (1986) 1238.
[14] A.M. Ishkhanyan, J. Phys. A **38** (2005) L491.
[15] T.A. Ishkhanyan, Y. Pashayan-Leroy, M.R. Gevorgyan, C. Leroy, A.M. Ishkhanyan, J. Contemp. Phys. (Arm. Ac. Sci.) **51** (2016) 229.
[16] P.P. Fiziev, J. Phys. A **43** (2010) 035203.
[17] A. Ishkhanyan and K.-A. Suominen, J. Phys. A **36** (2003) L81.
[18] M.N. Hounkonnou, A. Ronveaux, Appl. Math. Comput. **209** (2009) 421.
[19] C. Leroy and A.M. Ishkhanyan, Integral Transforms and Special Functions **26** (2015) 451.
[20] V.A. Shahnazaryan, T.A. Ishkhanyan, T.A. Shahverdyan, and A.M. Ishkhanyan, Armenian J. Phys. **5** (2012) 146.
[21] M. Hortacsu, Proc. 13th Regional Conference on Mathematical Physics, (World Scientific, Singapore, 2013) 23-39.
[22] A.Ya. Kazakov and S.Yu. Slavyanov, Theor. Math. Phys. **155** (2008) 722.
[23] S.Y. Slavyanov, Constr. Approx. **39** (2014) 75.
[24] K. Takemura, J. Math. Anal. Appl. **342** (2008) 52.
[25] S.Yu. Slavyanov, D.A. Satco, A.M. Ishkhanyan, T.A. Rotinyan, Theor. Math. Phys. **189** (2016) 1726.
[26] K. Iwasaki, H. Kimura, S. Shimomura, and M. Yoshida, From Gauss to Painlevé: A modern theory of special functions, Aspects of mathematics, v. 16 (Vieweg, Braunschweig, 1991).
[27] T.A. Ishkhanyan and A.M. Ishkhanyan, AIP Advances **4** (2014) 087132.